\documentclass{amsart}
\usepackage[left=3cm,right=3cm,top=2cm,bottom=2cm]{geometry}
\usepackage[fontsize=12pt]{fontsize}
\usepackage{graphicx}
\usepackage{cite}
\usepackage{mathtools}
\newtheorem{theorem}{Theorem}
\newtheorem*{theorem*}{Theorem}

\theoremstyle{definition}

\newtheorem{example}{Example}

\newtheorem{corollary}{Corollary} 
\newtheorem*{corollary*}{Corollary}
\theoremstyle{remark}

\numberwithin{equation}{section}

\begin{document}

\title{On the inverse M\"obius transformation and unrestricted partitions}


\author{Romulo L. Cruz-Simbron}
\address{Universidad Nacional de Ingenieria, Av. Tupac Amaru, 210, Rimac, Lima, Peru}
\address{Chemistry Department, University of Colorado Boulder, Boulder 80303, Colorado}
\curraddr{}
\email{romulo.cruz-simbron@colorado.edu}
\thanks{}



\date{\today}

\dedicatory{}

\begin{abstract}
\normalsize In this work we explore the relationship between the M\"obius inverse transformation of an arbitrary sequence and unrestricted partitions identities. 
\end{abstract}

\maketitle

\section{Introduction}
A partition  of a number $n$ is a sequence denoted by $\pi = (k_1, k_2, \cdots)$ comprising non-negative integers, with a finite subset of the $k_i$ being nonzero. $\pi$ is a partition of $n$ if and only if $n$ can be expressed as the sum $\sum jk_j$, where each positive integer $j$ appears with a frequency or multiplicity of $k_j$, without consideration of order. The integers $j$'s are called the parts of the partition. For example, the partition $\pi=1+2+3+2+1$  of 9 is denoted $\pi = (2,2,1,0,0,\cdots)$ with $k_1=2$, $k_2=2$ ,  $k_3=1$ and $k_4=k_5=\cdots = 0$.\\

Within this study, our focus lies on sums represented by $\sum_{\pi \in P(n)}f(\pi)$, wherein the summation encompasses all partitions of $n$, denoted by $P(n)$.\\

Four particular functions that we shall discuss are 
\begin{equation}
k(\pi) = \sum_j k_j
\end{equation}
 \noindent the total number of parts in $\pi$, and
 \begin{equation}
 Q(\pi) = \sum_{k_j>0} 1
 \end{equation}

\noindent the number of distinct parts in $\pi$. For example, the partition $\pi = (2,2,1)$ of 9 has $k(\pi)=5$ and $Q(\pi)=3$.\\

The other two function are the $s(\pi)$ and $l(\pi)$, the smallest and largest part in a partition $\pi$ of n.\\

A weighted partition identity is defined as a relation involving the sum $\sum_{\pi \in P(n)}f(\pi)$, where the function $f(\pi)$ incorporates functions of $k(\pi)$, $Q(\pi)$, $s(\pi)$, $l(\pi)$ or related variables. The investigation of weighted partition identities represents a relatively contemporary pursuit, laden with substantial implications yet to be fully elucidated \cite{berkovich2016new}. Beginning with the early studies of Uchimura \cite{uchimura1981identity} and Bressoud - Subbarao \cite{bressoud1984uchimura} on weighted partition identities concerning the sum of partitions into distinct parts, research has progressed systematically, notably under the direction of Alladi \cite{alladi1997partition, alladi1998partition, alladi2016partitions, alladi2023schmidt}, exploring weighted partition identities of unrestricted partitions. Subsequently, numerous results have emerged from the contributions of Andrews \cite{andrews2008number, andrews2013self}, Chen \cite{chen2007weighted}, Uncu \cite{berkovich2016new, berkovich2017variation,uncu2018weighted,bridges2022weighted}, Garvan \cite{garvan2017weighted}, Lovejoy \cite{choi2009partitions}, Merca \cite{merca2023additive,andrews2020number}, R. Fokkink - W. Fokkink - Wang \cite{bing1995relation}, and Dixi - Agarwal - Bhoria - Eyyunni - Maji \cite{agarwal2023bressoud,dixit2020partition,bhoria2022generalization,gupta2021sum}.\\

In a previous study \cite{cruz2023fine}, and using Fine's \cite{fine1988basic} tools, or a elementary proof \cite{simbron2023inverse},  we establish that the Bressoud - Subbarao type identities could be interpreted as the inverse Mobius transformation of an arbitrary sequence. Continue that methodology in this study we found an analog identity of weighted unrestricted partitions.  

\section{Arithmetical interpretation of weighted partitions identities related to the largest and smallest part} 

\begin{theorem}
Let $a_n$ an arbitrary sequence of real or imaginary terms, $n$ a positive integer number, and $t$ and $u$  real or imaginary numbers. Then
\begin{equation}
\sum_{n \geq 1} a_n \frac{tu \left( (1-u)tq^{n+1} \right)_{\infty}}{\left( t q^n \right)_{\infty} } q^n = \sum_{n\geq1} q^n \sum_{\pi \in P(n)}  t^{k(\pi)}u^{Q(\pi)}a_{s(\pi)} 
\end{equation}
where $s(\pi)$ is the smallest term in the partitions $\pi$
\end{theorem}

\begin{theorem}
Let $a_n$ an arbitrary sequence of real or imaginary terms, $n$ a positive integer number, and $t$ and $u$  real or imaginary numbers. Then
\begin{equation}
\sum_{n \geq 1} a_n \frac{tu \left( (1-u)tq \right)_{n-1}}{\left( t q \right)_{n} } q^n = \sum_{n\geq1} q^n \sum_{\pi \in P(n)}  t^{k(\pi)}u^{Q(\pi)}a_{l(\pi)} 
\end{equation}
where $l(\pi)$ is the largest term in the partitions $\pi$
\end{theorem}

\begin{theorem}
Let $a_n$ an arbitrary sequence of real or imaginary terms, $n$ a positive integer number, and $t$ and $u$  real or imaginary numbers. Then
\begin{equation}
\sum_{n \geq 1} a_n  \sum_{i \geq 0} \frac{tu \left( (1-u)tq^{i+1} \right)_{n-1}}{\left( t q^{i+1} \right)_{n} } q^{n+i} = \sum_{n\geq1} q^n \sum_{\pi \in P(n)}  t^{k(\pi)}u^{Q(\pi)} \sum_{j=1}^{s(\pi)}a_{l(\pi)-s(\pi)+j} 
\end{equation}
where $l(\pi)$ and $s(\pi)$ are the largest and smallest term in the partitions $\pi$, respectively. 
\end{theorem}

\begin{theorem}
Let $a_n$ and arbitrary sequence of real or imaginary terms , $t$ a real or imaginary number different of $1$ or $0$, and  $n$ a positive integer number, Then
\begin{equation}
\sum_{\pi \in P(n)}  (t)^{k(\pi)-1}\left(\frac{t-1}{t}\right)^{Q(\pi)-1} \sum_{j=1}^{s(\pi)}a_{l(\pi)-s(\pi)+j} = \sum_{d \vert n}  (t)^{d} \left(\frac{a_1}{t}+\frac{a_2}{t^2}+\frac{a_3}{t^3}+ \cdots +\frac{a_d}{t^d}\right)
\end{equation}
\label{MainTheorem}
\end{theorem}

\begin{corollary}
Let $a_n$ and arbitrary sequence of real or imaginary terms , $t$ a real or imaginary number different of $1$ or $0$, and  $n$ a positive integer number, Then
\begin{align}
\begin{split}
\sum_{\pi \in P(n)}  (t)^{k(\pi)-1}\left( \frac{t-1}{t} \right)^{Q(\pi)-1} &\sum_{j=1}^{s(\pi)} (t)^{l(\pi)-s(\pi)+j}a_{l(\pi)-s(\pi)+j} \\
&= \sum_{d \vert n}  (t)^{d} (a_1+a_2+a_3+ \cdots +a_d)
\end{split}
\end{align}
\label{corollary1}
\end{corollary}

\begin{corollary}
Let $a_n$ and arbitrary sequence of real or imaginary terms with $a_0=0$ and  $t$ a real or imaginary number different of $1$ or $0$,. 
\begin{align}
\begin{split}
\sum_{\pi \in P(n)}  (t)^{k(\pi)-1}\left(\frac{t-1}{t}\right)^{Q(\pi)-1} &\sum_{j=1}^{s(\pi)}  \left( a_{l(\pi)-s(\pi)+j} - t a_{l(\pi)-s(\pi)+j-1} \right)  = \sum_{d \vert n} a_d
\end{split}
\end{align}
\label{corollary2}
\end{corollary}
We should consider that $a_0=0$ can be defined for any sequence $a_1, a_2, \cdots$, simply by appending zero to it: $0, a_1, a_2, \cdots$. This Corollary show that the Mobius inverse transformation of any sequence $a_n$ with $a_0=0$ could be interpreted as weighted partition. We should notice that the left therm depend in $t$ and the right term not.

\begin{example}
If in Theorem \ref{MainTheorem} $a_n= n$ and $t=-1$, Then
\begin{align}
\begin{split}
\sum_{\pi \in P(n)}  (-1)^{k(\pi)-1}2^{Q(\pi)-1} &\sum_{j=1}^{s(\pi)} l(\pi)-s(\pi)+j \\
&= \sum_{d \vert n}  (-1)^{d} \left(\frac{1}{-1}+\frac{2}{(-1)^2}+\frac{3}{(-1)^3}+ \cdots +\frac{d}{(-1)^d}\right)\\
&= \sum_{d \vert n}  (-1)^{d} \left(-1+2-3+ \cdots +(-1)^dd\right)\\
&= \sum_{d \vert n} \left[  \frac{d}{2} + \left(\frac{1-(-1)^d}{2}\right)\frac{1}{2} \right]\\
&= \frac{\sigma_1(n) + \tau_{odd}(n)}{2} 
\end{split}
\end{align}
\end{example}

\begin{example}
If in Theorem \ref{MainTheorem} $a_n= n^2$ and $t=-1$, Then
\begin{align}
\begin{split}
\sum_{\pi \in P(n)}  (-1)^{k(\pi)-1}2^{Q(\pi)-1} &\sum_{j=1}^{s(\pi)} (l(\pi)-s(\pi)+j)^2 \\
&= \sum_{d \vert n}  (-1)^{d} \left(\frac{1^2}{-1}+\frac{2^2}{(-1)^2}+\frac{3^2}{(-1)^3}+ \cdots +\frac{d^2}{(-1)^d}\right)\\
&= \sum_{d \vert n}  (-1)^{d} \left(-1+4-9+ \cdots +(-1)^dd^2\right)\\
&= \sum_{d \vert n} \frac{d(d+1)}{2} \\
&= \frac{\sigma_2(n) + \sigma_{1}(n)}{2} 
\end{split}
\end{align}
\end{example}

\begin{example}
If in the Corollary \ref{corollary2} $a_{2k}=0, a_{2k+1}=1 $ (for $k \geq 0$) and $t=-1$, Then
\begin{align}
\begin{split}
&\sum_{\pi \in P(n)}  (-1)^{k(\pi)-1}2^{Q(\pi)-1} \sum_{j=1}^{s(\pi)} 1 = \sum_{d \vert n} \frac{1 -(-1)^d}{2}  \\
&\sum_{\pi \in P(n)}  (-1)^{k(\pi)-1}2^{Q(\pi)-1} s(\pi)= \tau_{odd}(n) \\
\end{split}
\end{align}
\end{example}
where $t_{odd}$ count the odd divisor of n.
\section{Proof of the Theorems and Corollaries}
\subsection*{Proof of Theorem 1}
\begin{proof}

Nathan Fine [Theorem 1, Chapter 2 in \cite{fine1988basic}] has establish that if $\psi_j(q) =  \sum_{ k \geq 0} C_j(k)q^k  \; (j = 1,2,3,\cdots)$. Then
\begin{equation}
\prod_{j\geq1} \psi_j(q^j) = \sum_{n\geq0} q^n \sum_{\pi \in P(n)} C_1(k_1)C_2(k_2)C_3(k_3) \cdots
\end{equation}

This Theorem allows to get straightforward arithmetical interpretations. In fact , the following product
\begin{align*}
(1+tuq+t^2uq^2+\cdots)(1+tuq^2+t^2uq^4+\cdots)(1+tuq^3+t^2uq^6+\cdots)\cdots\\
\end{align*}
has a $C_j(k) = t^k u $ ($k\geq 1$) and $C_j(0)=1$. Then, it means that $ \sum_{n\geq1} q^n \sum_{\pi \in P(n)}$ $ C_1(k_1)C_2(k_2)C_3(k_3)\cdots = $ $ \sum_{n\geq1} q^n \sum_{\pi \in P(n)}  t^{k(\pi)}u^{Q(\pi)}$, where $k(\pi)$ and $Q(\pi)$ represent the number of parts and the number of distinct parts of the partitions, respectively. Now, if in the preceding product we exclude the 1 from the first factor,
\begin{align*}
(0+tuq+t^2uq^2+\cdots)(1+tuq^2+t^2uq^4+\cdots)(1+tuq^3+t^2uq^6+\cdots)\cdots\\
\end{align*}
it has $C_j(k) = t^k u $ ($k\geq 1$) , $C_1(0)=0$ and $C_j(0)=1$ (for $j \geq 2$) . Then, it means that $ \sum_{n\geq1} q^n \sum_{\pi \in P(n)}$  $C_1(k_1)C_2(k_2)C_3(k_3)\cdots = \sum_{n\geq1} q^n \sum_{\pi \in P(n)}  t^{k(\pi)}u^{Q(\pi)}$ $ [k_1 \neq 0]$, where the symbol $[k_1 \neq 0]$ means that the smallest term must be $1$ in the partitions of n we consider for this sum. It should be noted that the summation begins from 1 instead of 0, as all the partitions must have at least the part 1 ($[k_1 \neq 0]$). Now, we can generalize this interpretation if we have the product

\begin{align*}
(0+tuq^m+t^2uq^{2m}+\cdots)&(1+tuq^{m+1}+t^2uq^{2(m+1)}+\cdots)\\
&\times(1+tuq^{m+2}+t^2uq^{2(m+2)}+\cdots)\cdots\\
\end{align*}
which can be rewrite as 
\begin{align*}
(1+0q^1+0q^2+\cdots)&(1+0q^2+0q^4+\cdots) \\
&\times(1+0q^3+0q^6+\cdots)(1+0q^4+0q^8+\cdots) \\
&\vdots\\
&\times(0+tuq^m+t^2uq^{2m}+\cdots)(1+tuq^{m+1}+t^2uq^{2(m+1)}+\cdots)
\end{align*}

\noindent it has $C_j(k) = 0$ ($k\geq 1$, $j= 1 \cdots (m-1)$)  , $C_j(0) = 1$ ( $j= 1 \cdots (m-1)$) , $C_j(k) = t^k u $ ($k\geq 1$,$j \geq m$ ) , $C_m(0)=0$ and $C_j(0)=1$ (for $j \geq (m+1)$). Using the Fine's Theorem 1, we got that this product is equal to $ \sum_{n\geq1} q^n \sum_{\pi \in P(n)}$  $C_1(k_1)C_2(k_2)C_3(k_3)\cdots = \sum_{n\geq1} q^n \sum_{\pi \in P(n)}  t^{k(\pi)}u^{Q(\pi)}$ $ [k_1 = 0 \wedge k_2 = 0 \cdots \wedge k_{m-1} = 0 $ $\wedge k_m \neq 0]$, where the symbol $[\quad ]$ represent the condition of this sum. This conditional could be interpreted as the sum go over all the partitions of n which smallest term is $m$.

Based on the previous result, if we have the following sum 

\begin{align}
\begin{split}
&a_1(0+tuq+t^2uq^2+\cdots)(1+tuq^2+t^2uq^4+\cdots)(1+tuq^3+t^2uq^6+\cdots)\cdots\\
&+a_2(0+tuq^2+t^2uq^4+\cdots)(1+tuq^3+t^2uq^6+\cdots)(1+tuq^4+t^2uq^8+\cdots)\cdots\\
&+a_3(0+tuq^3+t^2uq^6+\cdots)(1+tuq^4+t^2uq^8+\cdots)(1+tuq^5+t^2uq^{10}+\cdots)\cdots\\
&+a_4(0+tuq^4+t^2uq^8+\cdots)(1+tuq^5+t^2uq^{10}+\cdots)(1+tuq^6+t^2uq^{12}+\cdots)\cdots\\
&+a_5(0+tuq^5+t^2uq^{10}+\cdots)(1+tuq^6+t^2uq^{12}+\cdots)(1+tuq^7+t^2uq^{14}+\cdots)\cdots\\
&+a_6(0+tuq^6+t^2uq^{12}+\cdots)(1+tuq^7+t^2uq^{14}+\cdots)(1+tuq^8+t^2uq^{16}+\cdots)\cdots\\
\vdots
\end{split}
\label{sum_smallest}
\end{align}

it could be interpreted as 
\begin{align*}
\sum_{n\geq1} q^n \sum_{\pi \in P(n)} & t^{k(\pi)}u^{Q(\pi)}(a_1[k_1 \neq 0] + a_2[k_1 = 0 \wedge k_2 \neq 0] + \cdots ) \\
&= \sum_{n\geq1} q^n \sum_{\pi \in P(n)}  t^{k(\pi)}u^{Q(\pi)}a_{s(\pi)} 
\end{align*} 

This allow us to prove the theorem. \noindent In fact, the sum  in \ref{sum_smallest} is 
\begin{align*}
a_1 [ tuq &(1+tq + t^2q^2 + \cdots) ]\\
&\times[1+tuq^2(1+tq^2+t^2q^4+\cdots)]\\
&\times[1+tuq^3(1+tq^3+t^2q^6+\cdots)]\\
&\vdots\\
&+\\
a_2 [tuq^2 &(1+tq^2 + t^2q^4 + \cdots)]\\
&\times[1+tuq^3(1+tq^3+t^2q^6+\cdots)]\\
&\times[1+tuq^4(1+tq^4+t^2q^8+\cdots)]\\
&\vdots\\
&+\\
a_3 [tuq^3 &(1+tq^3 + t^2q^6 + \cdots) ]\\
&\times[1+tuq^4(1+tq^4+t^2q^8+\cdots)]\\
&\times[1+tuq^5(1+tq^4+t^2q^{10}+\cdots)]\\
&\vdots\\
&+\\
&\vdots
\end{align*}
\noindent which is equal to 
\begin{align*}
&a_1 \left( \frac{tuq}{1-tq}\right) \left( \frac{1-(1-u)tq^2}{1-tq^2}\right)\left( \frac{1-(1-u)tq^3}{1-tq^3}\right) \cdots\\
&+\\
&a_2 \left( \frac{tuq^2}{1-tq^2}\right) \left( \frac{1-(1-u)tq^3}{1-tq^3}\right)\left( \frac{1-(1-u)tq^4}{1-tq^4}\right) \cdots\\
&+\\
&a_3 \left( \frac{tuq^3}{1-tq^3}\right) \left( \frac{1-(1-u)tq^4}{1-tq^4}\right)\left( \frac{1-(1-u)tq^5}{1-tq^5}\right) \cdots\\
&+\\
&\vdots
\end{align*}
\noindent If we summarize this sum and represent the product as a geometric series, we obtain the desired result.

\end{proof}

\subsection*{Proof of Theorem 2}
\begin{proof}
Similarly to the previous proof, we start with the analysis of the product:

\begin{align*}
(1+tuq+&t^2uq^{2}+\cdots)(1+tuq^{2}+t^2uq^{4}+\cdots)\cdots\\
&\times(1+tuq^{m-1}+t^2uq^{2(m-1)}+\cdots)(0+tuq^m+t^2uq^{2m}+\cdots)
\end{align*}
which can be rewrite as 
\begin{align*}
(1+tuq+&t^2uq^{2}+\cdots)(1+tuq^{2}+t^2uq^{4}+\cdots)\cdots\\
&\times(1+tuq^{m-1}+t^2uq^{2(m-1)}+\cdots)(0+tuq^m+t^2uq^{2m}+\cdots)\\
&\times(1+0q^{m+1}+0q^{2(m+1)}+\cdots)(1+0q^{m+2}+0q^{2(m+2)}+\cdots)\cdots\\
\end{align*}

\noindent In this generation function, $C_j(0)= 1$ ( for $j = 1 \cdots m-1 $), $C_j(0)= 0$ ( for $j = m $), $C_j(0)= 1$ ( for $j \geq m+1 $) , $C_j(k) = t^ku$ (For $k\geq 1$, $j=1 \cdots m$) and $C_j(k) = 0$ (For $k\geq 1$, $j \geq m+1$). Using the Fine's Theorem 1, we got that this product is equal to $ \sum_{n\geq1} q^n \sum_{\pi \in P(n)}$  $C_1(k_1)C_2(k_2)C_3(k_3)\cdots = \sum_{n\geq1} q^n \sum_{\pi \in P(n)}  t^{k(\pi)}u^{Q(\pi)}$ $ [k_m \neq 0 \wedge k_{m+1} = 0 \wedge k_{m+2} = 0 $ $\wedge \cdots]$. This conditional could be interpreted as the sum go over all the partitions of n which largest term is $m$.

Based on the previous result, if we have the following sum 

\begin{align}
\begin{split}
&a_1(0+tuq+t^2uq^2+\cdots)\\
&+a_2(1+tuq+t^2uq^2+\cdots)(0+tuq^2+t^2uq^4+\cdots)\\
&+a_3(1+tuq+t^2uq^2+\cdots)(1+tuq^2+t^2uq^4+\cdots)(0+tuq^3+t^2uq^6+\cdots)\\
\vdots
\end{split}
\label{sum_largest}
\end{align}

it could be interpreted as 
\begin{align*}
\sum_{n\geq1} q^n \sum_{\pi \in P(n)} & t^{k(\pi)}u^{Q(\pi)}(a_1[k_1 \neq 0 \wedge k_{2} = 0 \wedge \cdots ] + a_2[k_2 \neq 0 \wedge k_{3} = 0 \wedge \cdots ] + \cdots ) \\
&= \sum_{n\geq1} q^n \sum_{\pi \in P(n)}  t^{k(\pi)}u^{Q(\pi)}a_{l(\pi)} 
\end{align*}

The sum  \ref{sum_largest} is 
\begin{align*}
a_1 [ tuq &(1+tq + t^2q^2 + \cdots) ]\\
&+\\
&a_2 [ 1+ tuq (1+tq + t^2q^2 + \cdots) ]\\
&\times [tuq^2 (tq^2 + t^2q^4 + \cdots)]\\
&+\\
&a_3 [ 1+ tuq (1+tq + t^2q^2 + \cdots) ]\\
&\times [1+tuq^2 (tq^2 + t^2q^4 + \cdots)]\\
&\times [tuq^3 (1+tq^3 + t^2q^6 + \cdots) ]\\
&+\\
&\vdots
\end{align*}
\noindent which is equal to 
\begin{align*}
&a_1 \left( \frac{tuq}{1-tq}\right)\\
&+\\
&a_2 \left( 1+ \frac{tuq}{1-tq}\right) \left( \frac{tuq^2}{1-tq^2}\right) \\
&+\\
&a_3 \left( 1+ \frac{tuq}{1-tq}\right) \left( 1+\frac{tuq^2}{1-tq^2}\right) \left( \frac{tuq^3}{1-tq^3}\right)\\
&+\\
&\vdots
\end{align*}
\noindent If we summarize this sum and represent the product as a geometric series, we obtain the desired result.

\end{proof}

\subsection*{Proof of Theorem 3}
\begin{proof}
Continuing with our methodology, we start with the product 
\begin{align*}
(1+tuq^{m-p+1}+&t^2uq^{2(m-p+1)}+\cdots)(1+tuq^{m-p+2}+t^2uq^{2(m-p+2)}+\cdots)\cdots\\
&\times(1+tuq^{m-1}+t^2uq^{2(m-1)}+\cdots)(0+tuq^m+t^2uq^{2m}+\cdots)
\end{align*}
This product is similar to the product we had explored for the largest-term-related partitions , but in this case, it is truncated into $p$ factor. We can rewrite this product as 
\begin{align*}
(1+&0q^1+0q^{2}+\cdots)(1+0q^{2}+0q^{4}+\cdots)\cdots\\
&\times(1+0q^{m-p-1}+0q^{m-p-1}+\cdots)(1+0q^{m-p}+0q^{m-p}+\cdots) \\
&\times(1+tuq^{m-p+1}+t^2uq^{2(m-p+1)}+\cdots)(1+tuq^{m-p+2}+t^2uq^{2(m-p+2)}+\cdots)\cdots \\
&\times(1+tuq^{m-1}+t^2uq^{2(m-1)}+\cdots)(0+tuq^m+t^2uq^{2m}+\cdots)\\
&\times(1+0q^{m+1}+0q^{2(m+1)}+\cdots)(1+0q^{m+2}+0q^{2(m+2)}+\cdots)\cdots\\
\end{align*}

\noindent In this generation function, $C_j(0)= 1$ ( for $j = 1, \cdots, m-1 $), $C_j(0)= 0$ ( for $j = m $), $C_j(0)= 1$ ( for $j \geq m+1 $) , $C_j(k) = 0$ (For $k\geq 1$, $j=1 \cdots m-p$), $C_j(k) = t^ku$ (For $k\geq 1$, $j=m-p+1, \cdots , m$)  , and $C_j(k) = 0$ (For $k\geq 1$, $j \geq m+1$). Using the Fine's Theorem 1, we got that this product is equal to $ \sum_{n\geq1} q^n \sum_{\pi \in P(n)}$  $C_1(k_1)C_2(k_2)C_3(k_3)\cdots = \sum_{n\geq1} q^n \sum_{\pi \in P(n)}  t^{k(\pi)}u^{Q(\pi)}$ $ [\cdots \wedge $ $k_{m-p} = 0  $ $\wedge k_{m} \neq 0 \wedge k_{m+1} = 0 $ $\wedge \cdots]$. This conditional could be interpreted as the sum go over all the partitions of n which largest term must be $m$ and its smallest term could be $m-p+1$ or larger, but not smaller.

Based in the previous result, if we have the following sum 

\begin{align}
\begin{split}
&a_1(0+tuq+t^2uq^2+\cdots)\\
&+\\
&+a_2(1+tuq+t^2uq^2+\cdots)(0+tuq^2+t^2uq^4+\cdots)\\
&+a_1(0+tuq^2+t^2uq^4+\cdots)\\
&+\\
&+a_3(1+tuq+t^2uq^2+\cdots)(1+tuq^2+t^2uq^4+\cdots)(0+tuq^3+t^2uq^6+\cdots)\\
&+a_2(1+tuq^2+t^2uq^4+\cdots)(0+tuq^3+t^2uq^6+\cdots)\\
&+a_1(0+tuq^3+t^2uq^6+\cdots)\\
&+\\
&+a_4(1+tuq+t^2uq^2+\cdots)(1+tuq^2+t^2uq^4+\cdots)(1+tuq^3+t^2uq^6+\cdots)\\
&\times(0+tuq^4+t^2uq^8+\cdots)\\
&+a_3(1+tuq^2+t^2uq^4+\cdots)(1+tuq^3+t^2uq^6+\cdots)(0+tuq^4+t^2uq^8+\cdots)\\
&+a_2(1+tuq^3+t^2uq^6+\cdots)(0+tuq^4+t^2uq^8+\cdots)\\
&+a_1(0+tuq^4+t^2uq^8+\cdots)\\
&+\\
\vdots
\end{split}
\label{sum_both}
\end{align}

it could be represented as
\begin{align*}
\sum_{n\geq1} q^n \sum_{\pi \in P(n)}&  t^{k(\pi)}u^{Q(\pi)} (a_1[k_{1} \neq 0 \wedge k_{2} = 0  \wedge \cdots] \\
&+\\
& a_2[k_{2} \neq 0 \wedge k_{3} = 0  \wedge \cdots] +\\
& a_1[k_1 = 0 \wedge k_{2} \neq 0 \wedge k_{3} = 0  \wedge \cdots] \\
&+\\
& a_3[k_{3} \neq 0 \wedge k_{4} = 0  \wedge \cdots] +\\
& a_2[k_{1} = 0  \wedge k_{3} \neq 0 \wedge k_{4} = 0  \wedge \cdots] +\\
& a_1[k_{1} = 0  \wedge k_{2} = 0  \wedge k_{3} \neq 0 \wedge k_{4} = 0  \wedge \cdots] \\
&+\\
& a_4[k_{4} \neq 0 \wedge k_{5} = 0  \wedge \cdots] +\\
& a_3[k_{1} = 0  \wedge k_{4} \neq 0 \wedge k_{5} = 0  \wedge \cdots] +\\
& a_2[k_{1} = 0  \wedge k_{2} = 0  \wedge k_{4} \neq 0 \wedge k_{5} = 0  \wedge \cdots] +\\
& a_1[k_{1} = 0  \wedge k_{2} = 0  \wedge k_{3} = 0  \wedge k_{4} \neq 0 \wedge k_{5} = 0  \wedge \cdots] \\
&\cdots)\\
\end{align*} 

\noindent which could be rewrite and interpreted as 

\begin{align*}
\sum_{n\geq1} q^n \sum_{\pi \in P(n)}&  t^{k(\pi)}u^{Q(\pi)} (a_1[k_{1} \neq 0 \wedge k_{2} = 0  \wedge \cdots]+ \\
& a_1[k_1 = 0 \wedge k_{2} \neq 0 \wedge k_{3} = 0  \wedge \cdots]+ \\
& a_1[k_{1} = 0  \wedge k_{2} = 0  \wedge k_{3} \neq 0 \wedge k_{4} = 0  \wedge \cdots]+ \\
& a_1[k_{1} = 0  \wedge k_{2} = 0  \wedge k_{3} = 0  \wedge k_{4} \neq 0 \wedge k_{5} = 0  \wedge \cdots] +\\
&\vdots\\
&+\\
& a_2[k_{2} \neq 0 \wedge k_{3} = 0  \wedge \cdots] +\\
& a_2[k_{1} = 0  \wedge k_{3} \neq 0 \wedge k_{4} = 0  \wedge \cdots] +\\
& a_2[k_{1} = 0  \wedge k_{2} = 0  \wedge k_{4} \neq 0 \wedge k_{5} = 0  \wedge \cdots] +\\
&\vdots\\
&+\\
& a_3[k_{3} \neq 0 \wedge k_{4} = 0  \wedge \cdots] +\\
& a_3[k_{1} = 0  \wedge k_{4} \neq 0 \wedge k_{5} = 0  \wedge \cdots] +\\
&\vdots\\
&+\\
& a_4[k_{4} \neq 0 \wedge k_{5} = 0  \wedge \cdots] +\\
&\cdots)\\
&= \sum_{n\geq1} q^n \sum_{\pi \in P(n)}  t^{k(\pi)}u^{Q(\pi)} \sum_{j \geq 1}^{s(\pi)}a_{l(\pi)-s(\pi)+j} 
\end{align*} 

The sum  \ref{sum_both} is also equal to 
\begin{align*}
&a_1 \left[ \left( \frac{tuq}{1-tq} \right) + \left( \frac{tuq^2}{1-tq^2} \right) + \cdots \right]\\
&+\\
&a_2 \left[ \left( 1+ \frac{tuq}{1-tq}\right) \left( \frac{tuq^2}{1-tq^2}\right) + \left( 1+ \frac{tuq^2}{1-tq^2}\right) \left( \frac{tuq^3}{1-tq^3}\right) + \cdots\right] \\
&+\\
&a_3 \left[ \left( 1+ \frac{tuq}{1-tq}\right) \left( 1+\frac{tuq^2}{1-tq^2}\right) \left( \frac{tuq^3}{1-tq^3}\right) + \left( 1+ \frac{tuq^2}{1-tq^2}\right) \left( 1+\frac{tuq^3}{1-tq^3}\right) \left( \frac{tuq^4}{1-tq^4}\right) + \cdots\right]\\
&+\\
&\vdots
\end{align*}

\noindent If we summarize this sum and represent the product as a geometric series, we obtain the desired result.
\end{proof}

\subsection*{Proof of Theorem 4}
\begin{proof}
In fact, the left term in the Theorem 3 is
\begin{align*}
\sum_{n \geq 1} a_n  &\sum_{i \geq 0} \frac{tu \left( (1-u)tq^{i+1} \right)_{n-1}}{\left( t q^{i+1} \right)_{n} } q^{n+i}  = tu \sum_{n \geq 1} a_n  \sum_{i \geq 0} \frac{ (tq)_i \left( (1-u)tq \right)_{n+i-1}}{ (tq)_{n+i}\left( (1-u)tq \right)_{i} } q^{n+i} \\
& = tu \sum_{n \geq 1} a_n  q^n \frac{ \left( (1-u)tq\right)_{n-1} }{(tq)_n}  \sum_{i \geq 0} \frac{ (tq)_i \left( (1-u)tq^n \right)_{i}}{ (tq^{n+1})_{i}\left( (1-u)tq \right)_{i} } q^{i} \\
\end{align*}

If we set $t= t$ , $u =(\frac{t-1}{t})$, and we use the \textit{Heine’s transformation} as Andrews in \cite{andrews2015partitions}, 

\begin{equation}
\sum_{i \geq 0} \frac{(a)_i (b)_i}{(q)_i(c)_i}z^i = \frac{\left( \frac{c}{b} \right)_{\infty} (bz)_{\infty}}{(c)_{\infty} (z)_{\infty}} \sum_{j \geq 0} \frac{\left( \frac{abz}{c} \right)_{j} (b)_j \left( \frac{c}{b} \right)^j}{(q)_j(bz)_j}
\end{equation}
\noindent we got, 
\begin{align*}
(tu) \sum_{n \geq 1} a_n  & q^n \frac{ \left(q\right)_{n-1} }{(tq)_n}  \sum_{i \geq 0} \frac{ (tq)_i \left( q^n \right)_{i}}{ (tq^{n+1})_{i}\left( q \right)_{i} } q^{i} =\\
& (tu) \sum_{n \geq 1} a_n  q^n \frac{ \left(q\right)_{n-1} }{(tq)_n}  \frac{(tq)_{\infty}(q^{n+1})_{\infty}}{(tq^{n+1})_{\infty}(q)_{\infty}} \sum_{j \geq 0} \frac{ (q)_j \left( q^n \right)_{j}}{ (q)_{j}\left( q^{n+1} \right)_{j} }(tq)^j\\
& (tu) \sum_{n \geq 1} a_n  q^n \frac{ \left(q\right)_{n-1} }{(tq)_n}  \frac{(tq)_{n}}{(q)_{n}} \sum_{j \geq 0} \frac{(1-q^n)}{(1-q^{n+j})} (tq)^{j}\\
& (tu) \sum_{n \geq 1} a_n  \sum_{j \geq 0} \frac{(t)^{j}(q)^{n+j}}{(1-q^{n+j})} \\
&= (tu) a_1 \left[ \frac{(t)^0q}{1-q} + \frac{(t)^1q^2}{1-q^2} + \frac{(t)^2q^3}{1-q^3} \cdots  \right]+\\
&(tu) a_2 \left[ \frac{(t)^0q^2}{1-q^2} + \frac{(t)^1q^3}{1-q^3} + \frac{(t)^2q^4}{1-q^4} \cdots  \right]+\\
&(tu) a_3 \left[ \frac{(t)^0q^3}{1-q^3} + \frac{(t)^1q^4}{1-q^4} + \frac{(t)^2q^5}{1-q^5} \cdots  \right]+\\
&\vdots\\
&=(tu)\sum_{k=1}\left(\sum_{i\geq1}^{k}(t)^{k-i} a_i \right) \frac{q^k}{1-q^k}\\
&=(tu)\sum_{k=1} q^k \sum_{d \vert k} \left( \sum_{i\geq1}^{d}(t)^{d-i} a_i \right) \\
&=(tu)\sum_{k=1} q^k \sum_{d \vert k}  (t)^{d}\sum_{i\geq1}^{d}(t)^{-i} a_i \\
&=(tu)\sum_{k=1} q^k \sum_{d \vert k}  (t)^{d} \left(\frac{a_1}{t}+\frac{a_2}{t^2}+\frac{a_3}{t^3}+ \cdots +\frac{a_d}{t^d}\right) \\
\end{align*}
which conclude the proof.
\end{proof}

\subsection*{Proof of Corollary 1}
\begin{proof}
In effect, If in the Theorem \ref{MainTheorem} we consider an auxiliary sequence $b_n$ defined by 
\begin{equation}
a_n = (t)^nb_n
\end{equation}
replacing $a_n$ in the Theorem \ref{MainTheorem} we got 
\begin{equation}
\sum_{\pi \in P(n)}  (t)^{k(\pi)-1}\left( \frac{t-1}{t}\right)^{Q(\pi)-1} \sum_{j=1}^{s(\pi)}(t)^{l(\pi)-s(\pi)+j} b_{l(\pi)-s(\pi)+j} = \sum_{d \vert n}  (t)^{d} (b_1+b_2+b_3+ \cdots +b_d)
\end{equation}
which conclude the proof. 
\end{proof}
\subsection*{Proof of Corollary 2}
\begin{proof}
If we define the auxiliary sequence $b_n$ defined as 
\begin{equation}
b_n = a_1 + a_2 + \cdots + a_n
\end{equation}
with $b_0=0$, then
\begin{equation}
a_n = b_n-b_{n-1}
\end{equation}
for $n\geq 1$

In the Corollary \ref{corollary1}, we conclude then, 
\begin{align}
\begin{split}
&\sum_{\pi \in P(n)}  (t)^{k(\pi)-1}\left(\frac{t-1}{t}\right)^{Q(\pi)-1} \sum_{j=1}^{s(\pi)} (t)^{l(\pi)-s(\pi)+j} \left[ b_{l(\pi)-s(\pi)+j} - b_{l(\pi)-s(\pi)+j-1} \right] \\
&= \sum_{d \vert n}  (t)^{d} b_d
\end{split}
\end{align}
If we then define an auxiliary sequence $c_n$ as 
\begin{align*}
b_n = (t)^{-n}c_n
\end{align*}
Then, the interior sum is 
\begin{align}
\begin{split}
&\sum_{j=1}^{s(\pi)} (t)^{l(\pi)-s(\pi)+j} \left[ (t)^{-(l(\pi)-s(\pi)+j)}c_{l(\pi)-s(\pi)+j}  - (t)^{-(l(\pi)-s(\pi)+j-1)}c_{l(\pi)-s(\pi)+j-1} \right] \\
&\sum_{j=1}^{s(\pi)} \left[ c_{l(\pi)-s(\pi)+j} - t c_{l(\pi)-s(\pi)+j-1} \right] \\
\end{split}
\end{align}
Finally, 
\begin{equation}
\begin{split}
&\sum_{\pi \in P(n)}  (t)^{k(\pi)-1}\left(\frac{t-1}{t} \right)^{Q(\pi)-1} \sum_{j=1}^{s(\pi)} \left[ c_{l(\pi)-s(\pi)+j} - t c_{l(\pi)-s(\pi)+j-1} \right] \\
&= \sum_{d \vert n}  c_d
 \end{split}
\end{equation}
which conclude the proof. 
\end{proof}
\section*{Acknowledgments}
RC extends gratitude to Ruth Quispe Pilco, Julia Simbron Pacheco, and Antonio Padua Cruz Claros for their support during the development of this work. Additionally, RC expresses thanks to the PhD Program at CU Boulder for the scholarship supporting their PhD studies.

\bibliographystyle{unsrt}
\bibliography{my_bibtex}
\end{document}